\documentclass[12pt]{article}
\usepackage{amssymb,latexsym,amsmath}
\begin{document}
\title{An analogue of a  conjecture of Mazur: a question in 
Diophantine approximation on tori}
\author{Dipendra Prasad}
\maketitle
\newtheorem{thm}{Theorem}
\newtheorem{cor}{Corollary}
\newtheorem{prop}{Proposition}
\newtheorem{lemma}{Lemma}
\newtheorem{defn}{Definition}
\newtheorem{conj}{Conjecture}
\baselineskip=16pt
{\sf 

\hfill{\small To Professor Shalika, with admiration}

\begin{abstract}

B. Mazur has considered the question of density in the
Euclidean topology of the set of ${\Bbb Q}$-rational
points on a  variety $X$ defined over ${\Bbb Q}$, in particular
for Abelian varieties.  In this paper we consider the question
of closures of the image of finitely generated 
subgroups of $T({\Bbb Q})$ in $\Gamma \backslash T({\Bbb R})$
where $T$ is a torus defined over ${\Bbb Q}$, $\Gamma$ an arithmetic subgroup
such that $\Gamma \backslash T({\Bbb R})$ is compact. Assuming 
Schanuel's conjecture, we prove that the 
closures correspond to {\it algebraic} sub-tori of $T$.

\end{abstract}

Let $V$ be a smooth algebraic variety over ${\Bbb Q}$. The set $V({\Bbb R})$ 
acquires a topological structure from the Euclidean topology of ${\Bbb R}$.
It is known that $V({\Bbb  R})$ 
has finitely many connected components.
If $V({\Bbb  Q})$ is Zariski dense in $V$, it was conjectured by  B. Mazur, 
cf. [M1] and [M2], that the closure of
$V({\Bbb Q})$ in $V({\Bbb R})$ is a finite union
of connected components of $V({\Bbb  R})$. 
This conjecture was shown to be false in this generality
by Colliot-Th\'el\`ene, Skorobogatov, and Swinnerton-Dyer for an 
elliptic surface  who have proposed a slightly reformulated  conjecture, 
cf. [CSS]. However, the
present evidence seems to suggest that the following special case of 
Mazur's conjecture is true.

\begin{conj}(Mazur's conjecture for Abelian varieties): Let $A$ be  an 
abelian variety over ${\Bbb Q}$, and $G$ a subgroup of $A({\Bbb Q})$. 
Then the closure of $G$  in the 
Euclidean topology of $A({\Bbb R})$ contains
$B({\Bbb R})^0$ 
as a subgroup of finite index for a certain abelian subvariety
$B$ defined over ${\Bbb Q}$.
\end{conj}

The following theorem of M.Waldschmidt [W1] is the best result known towards
Mazur's conjecture.

\begin{thm}(Waldschmidt) If $A$ is a simple abelian variety over ${\Bbb Q}$
of dimension $d$ and if rank of $A({\Bbb Q})$ is $\geq d^2-d+1$, then 
the closure of $A({\Bbb Q}) $ in the Euclidean topology contains
$A({\Bbb R})^0$. 
\end{thm}

\section{The conjecture about tori}

In this section we propose the following analogue of Mazur's
conjecture for tori. 
We begin by recalling certain standard definitions.

Let $S$ be a torus defined over ${\Bbb Q}$, i.e., let $S$ be a commutative
linear algebraic group over ${\Bbb Q}$ which becomes isomorphic to
${\Bbb G}^n_m$ over the algebraic closure $\overline{\Bbb  Q}$ of
${\Bbb Q}$ for a  certain integer $n \geq 0$. The torus $S$ is called
isotropic over ${\Bbb   Q}$ if there exists a homomorphism of algebraic groups
$S \rightarrow {\Bbb G}_m$ defined over ${\Bbb Q}$. If $S$ is not
isotropic over ${\Bbb Q}$, it is called anisotropic. 

Given a linear algebraic group over ${\Bbb Q}$ such as $S$, it makes sense 
to talk of arithmetic subgroups  ${\Gamma}$ of $S({\Bbb Q})$. Any two
arithmetic subgroups $\Gamma_1$ and $\Gamma_2$ are commensurable, i.e.,
$\Gamma_1 \cap \Gamma_2$ is of finite index in both $\Gamma_1$ and
$\Gamma_2$. It is a consequence of 
Dirichlet unit theorem (or the general theorem due to Borel and
Harish-Chandra) that if $S$ is an anisotropic torus over ${\Bbb Q}$,
then $\Gamma\backslash S({\Bbb R})$ is a compact abelian group
for $\Gamma$ any arithmetic subgroup of $S({\Bbb Q})$. 
The connected
component of identity of 
$\Gamma\backslash S({\Bbb R})$ is a torus in the usual sense
of the word, i.e., a topological group isomorphic to $(S^1)^n$.

Unlike most other algebraic groups, tori have the property that there
is a unique {\it maximal} arithmetic group. This unique maximal
arithmetic subgroup is the subgroup $\Gamma$ 
of $T({\Bbb Q})$ defined as follows: 
$$\Gamma = 
\left \{ \gamma \in T({\Bbb Q}) \left | 
\begin{array}{l}\chi(\gamma) 
= {\rm~a~unit~in~the~ring~of~integers~of~} \bar{\Bbb Q}^*  \\
{\rm~for~all~characters~} \chi:
 T \rightarrow {\Bbb G}_m  
\end{array} \right. \right \}.$$

We now make the following analogue of Mazur's conjecture for tori.

\begin{conj}

Let $S$ be an algebraic anisotropic torus defined over ${\Bbb Q}$. 
Let $F$ be a finitely generated subgroup of $S({\Bbb Q})$. 
For $\Gamma$ any arithmetic subgroup of 
$S({\Bbb Q})$, the
connected component of identity of the closure of the image of $F$ in
$\Gamma \backslash S({\Bbb R})$ 
equals connected component of identity of
$\Gamma_T \backslash T({\Bbb R})$ for a certain subtorus $T$ of $S$
defined over ${\Bbb   Q}$ with $\Gamma_T = \Gamma \cap T({\Bbb  Q})$. 
Equivalently,
the closure of the image of $F$ is dense in the Euclidean topology in
the identity component of $\Gamma \backslash S({\Bbb R})$ if and only if 
any subgroup $G$ of $S({\Bbb Q})$ surjecting on the image of $F$ in
$\Gamma \backslash S({\Bbb R})$ is Zariski dense in $S$.
\end{conj}

\noindent{\bf Remark 1:} To any finitely generated subgroup $F$ 
of $S({\Bbb Q})$, there is a natural subtorus $T$ of $S$ 
over ${\Bbb Q}$ defined as the  (connected component of the identity of the) 
kernel of the group of characters
$$X_F = \{ \chi: S \rightarrow {\Bbb G}_m | \chi(a)= {\rm~a
~unit~in~the~ring~of~integers~of~}
\bar{\Bbb Q}^* {\rm~for~all~}a \in F \}.$$
By embedding an anisotropic torus in a product of norm 1 tori (as in Lemma 
3 below), it is easy to see that the torus $T$ defined here is the same as 
that which appears in the conjecture above.
We are sloppy here, as well as other places
in the paper, about making such assertions only up to connected components.

\vspace{4mm}

\noindent{\bf Remark 2:} Let the finitely generated group $F$ be generated
by $\{a_1,a_2,\cdots, a_n\}$. Since product of compact subgroups is compact,
it is clear that the closure of the image of $F$ in  
$\Gamma \backslash S({\Bbb R})$, $\Gamma$ an arithmetic subgroup of
$S({\Bbb R})$, is the product of the closures of cyclic groups generated
by $a_i$. {\it Hence it suffices to prove the conjecture for cyclic groups 
$F$}.  By embedding an anisotropic torus in a product of norm 1 tori 
(as in Lemma 3 below), we are  further reduced to proving the conjecture 
for a cyclic subgroup in a product of norm 1 tori. 

\vspace{4mm}

\noindent{\bf Remark 3:} The conjecture above can also be formulated
as the closure of a finitely generated subgroup in $S({\Bbb Q})$ containing
an arithmetic subgroup in $S({\Bbb R})$. We note that since an arithmetic
group in an anisotropic torus over ${\Bbb Q}$ which is split over ${\Bbb R}$
is Zariski dense, so will this finitely
generated subgroup. Thus in this sense, our conjecture has  a different
flavor than Mazur's conjecture. We note that there is also a conjecture,
different from the one formulated here, due
to Waldschmidt [W3, Conjecture 3.5 of Chapter 3] 
about the closure in the Euclidean 
topology of a finitely generated subgroup of $S({\Bbb Q})$.

\vspace{4mm}

\noindent{\bf Remark 4:} In this paper we will often 
be using without explicit mention the trivial remark that  
the connected component of identity 
of the closure in either Euclidean or Zariski topology of a finitely
generated subgroup $F$  (of a topological group or an Algebraic group) 
or a subgroup $G$ of $F$ of finite index is the same.

\vspace{4mm}

\noindent{\bf Example 1:} The simplest case of our conjecture is when
$T$ is a torus over ${\Bbb Q}$ such that $T({\Bbb R})$ itself is compact.
In this case we can take $T({\Bbb Z})$ to be the trivial subgroup 
of $T({\Bbb Q})$. We will thus be comparing the closures in Euclidean topology and
Zariski topology of a finitely generated subgroup of $T({\Bbb Q})$. That the
two closures are the same follows from the easily proven assertion that 
any continuous homomorphism from the compact group $T({\Bbb R})$ to
$S^1$ is the restriction to $T({\Bbb R})$ of an algebraic character
(defined over ${\Bbb C}$) from $T({\Bbb C})$ to ${\Bbb C}^*$.

\vspace{4mm}
 
\noindent{\bf Example 2:} 
Let $k$ be a totally real cubic extension
of ${\Bbb Q}$. Let $T = k^1$ denote the group of Norm 1 elements of
$k$. So $T({\Bbb Q})= k^1$, and the group of units of (the ring of integers 
of) $k$ of norm 1 
can be taken to be an arithmetic subgroup of $T({\Bbb Q})$. We have,
\begin{eqnarray*}
T({\Bbb R}) & = & \left[ k \otimes {\Bbb R} \right ]^1 \\
            & = & \left[ {\Bbb R} \oplus {\Bbb R} \oplus {\Bbb R} \right]^1 \\
            & \cong & {\Bbb R}^* \times {\Bbb R}^*. 
\end{eqnarray*}

We note that $(x,y) \rightarrow (\log x, \log y)$ gives an isomorphism
of the product of two copies of positive reals (under multiplication) with 
real numbers (under addition). We have thus a homomorphism from
$T({\Bbb  R})$ to ${\Bbb R} \times {\Bbb R}$ such that the image of 
the group of units of norm 1 is a discrete cocompact subgroup with quotient
$S^1 \times S^1$.  

It is easy to see that the torus $T$ of norm 1 elements of a cubic 
field has no nontrivial
subtorus defined over ${\Bbb Q}$. Our conjecture  in this
case will therefore say that any element of $T({\Bbb Q})$ 
which is not of finite order
will generate a dense subgroup of $T({\Bbb Z})\backslash T({\Bbb  R})$.
We make this  more concrete.

Observe that if $\Lambda \subset {\Bbb R}^2$ is a lattice in ${\Bbb R}^2$
and $v$ is a vector in ${\Bbb  R}^2$, then integral multiples of
$v$ is dense in $\Lambda \backslash {\Bbb R}^2$ if and only if there does 
not exist a nonzero  integer $q$, and a lattice point $\lambda \in \Lambda$ 
such that $qv+\lambda$ 
is real multiple of a vector in $\Lambda$.  Suppose if possible,
$qv+\lambda = r \lambda_1$ for an integer $q$, and a real number $r$.  
By looking at the co-ordinates of the vectors on the two sides of the
equality,  it is easily seen that it  suffices to prove
that if $(\log|x_1|, \log|x_2|)$ is a real multiple of 
$(\log|\epsilon_1|, \log|\epsilon_2|)$  
then it is a rational multiple of 
$(\log|\epsilon_1|, \log|\epsilon_2|)$.
Here $x_1$ and $x_2$ are the images of an
element (corresponding to $qv+r$) in $k^1$ under two fixed 
embeddings into the reals,   and
$\epsilon_1$ and $\epsilon_2$ are the images of a unit
element in $k^1$ under the same embeddings into the reals.

This will be a consequence of the following  well-known conjecture, cf. [W2].

\begin{conj}(4 Exponential Conjecture due to Schneider, Lang and Ramachandra)
Let $M$ be a $2 \times 2$ matrix consisting of logarithm of algebraic
numbers. Assume that the rows of the matrix are linearly independent
over ${\Bbb  Q}$,  and also that the columns of
 the matrix are linearly independent over ${\Bbb  Q}$, then the determinant 
of $M$ is non-zero.
\end{conj}

The following theorem is a step towards the proof of the 4 exponential
conjecture. 

\begin{thm}(Lang and Ramachandra) Let $M$ be a $2 \times 3$ matrix consisting of logarithm of 
algebraic numbers. Assume that the rows of the matrix are linearly independent
over ${\Bbb  Q}$,  and also that the columns of
 the matrix are linearly independent over ${\Bbb  Q}$. Then the rank 
of $M$ is 2.
\end{thm}

\section{Some lemmas about Tori} 

In this section we collect together some elementary lemmas about tori. We will
be considering closures of finitely generated subgroups in the 
Euclidean and Zariski topologies.

\begin{lemma} (a) For a discrete subgroup $\Lambda \subset {\Bbb R}^n$ with 
${\Bbb R}^n/\Lambda$ compact, the integral multiples of a point $x \in
{\Bbb R}^n$ are dense inside ${\Bbb R}^n/\Lambda$ if and only if no
non-trivial continuous homomorphism of ${\Bbb R}^n/\Lambda$ to $S^1$ 
takes $x$ to the identity element of $S^1$. 

(b) The integral multiples of $x \in
{\Bbb R}^n$ are dense inside ${\Bbb R}^n/\Lambda$ if and only if $rx+\lambda$ 
does not belong to $\Lambda_1\otimes_{\Bbb Z}{\Bbb R}$ for any subgroup
$\Lambda_1$ of $\Lambda$ with rank$_{\Bbb Z}\Lambda_1 < n$, any nonzero integer
$r$, and any element $\lambda \in \Lambda$.

(c) The integral multiples of $x =(x_1,\cdots,x_n) \in
{\Bbb R}^n$ are dense inside ${\Bbb R}^n/\Lambda$ if and only if for 
 $\ell_1=(\ell_{11},\cdots,\ell_{1n}),\cdots, \ell_{n-1}=(\ell_{n-1,1},\cdots,
\ell_{n-1,n})$, belonging to $\Lambda $ and generating a rank $(n-1)$ subgroup 
of $\Lambda$, the matrix 
$$
 \left(
\begin{array}{cccccc}

rx_1+\lambda_1 & rx_2+\lambda_2 & & \cdot   & \cdot  &rx_n+\lambda_n \\
  \ell_{11} & \ell_{12}&  & \cdot  &  &\ell_{1n} \\
\cdot  & \cdot & \cdot  & \cdot &  & \cdot \\
 \cdot &  \cdot &  & \cdot &  &\cdot \\
\cdot & \cdot & \cdot & & \cdot & \cdot\\
\ell_{n-1,1} &\ell_{n-1,2}  & \cdot &  &\cdot  & \ell_{n-1,n}

\end{array}
\right)
 $$
is non-singular, i.e., the determinant is non-zero, for any nonzero
integer $r$, and any $\lambda=(\lambda_1,\cdots,\lambda_n) \in \Lambda$.

\end{lemma}

Now we have a lemma about density of the abelian group generated by a point
on a torus in the Zariski topology. In this lemma for a number field
$K$, we will be looking at the torus $T= R_{K/{\Bbb Q}}{\Bbb G}_m$ defined 
over ${\Bbb Q}$ to be the algebraic group whose group of rational 
points over any ${\Bbb Q}$-algebra $A$ is
$T(A) = (K\otimes_{\Bbb Q}A)^*$; in particular $T({\Bbb Q})=K^*$. 

\begin{lemma} (a) An element $x \in T({\Bbb Q}) = K^*$ generates $K$ 
(i.e., $K$ is the smallest field extension of ${\Bbb Q}$ containing $x$)
if and only if all its conjugates (i.e., the image of $x$ under all the 
distinct embeddings of $K$ in ${\Bbb C}$) $\{x_1,\cdots ,x_n\}$ are distinct. 

(b) An element $x \in T({\Bbb Q})= K^*$ lies in no proper algebraic subgroup 
defined over ${\Bbb Q}$ if and only if the abelian subgroup generated by
$\{x_1,\cdots,x_n\}$ is free abelian of rank $n$.

\end{lemma}

\noindent{\bf Proof :} Part $(a)$ is clear. For part $(b)$ note that 
any algebraic character of $T({\Bbb Q}) =K^*$ 
(defined over the algebraic closure) is defined by $z \rightarrow 
z_1^{d_1}\cdot z_2^{d_2}\cdots z_n^{d_n}$, where $z_i$ denotes the image of 
$z$ under the various embeddings of $K$ into ${\Bbb C}$. Since for an element
belonging to a proper algebraic subgroup of $T$, there is a character of $T$ 
trivial on that element, therefore if $x$ belongs to a proper algebraic
subgroup, the subgroup generated by $\{x_1,\cdots,x_n\}$ will not be
free. 

Conversely, if the subgroup generated by $\{x_1,\cdots,x_n \}$ is not 
free abelian,
$x$ belongs to the kernel of a non-trivial character $\chi$ of $T$ defined
over $\overline{\Bbb Q}$. Hence $x \in T({\Bbb Q})$ lies in 
$S(\overline {\Bbb Q})$ 
for an algebraic subgroup  $S$ of $T$. By Galois conjugation, 
$x  \in S^{\sigma}(\overline {\Bbb Q})$ for all Galois conjugates of $S$. 
Hence $x  \in \cap (S^{\sigma})(\overline {\Bbb Q})$.  However,
$ A = \cap (S^{\sigma})$ is an algebraic group defined over ${\Bbb Q}$. Hence
$x$  belongs to $ A({\Bbb Q})$  for $A$ a proper algebraic subgroup of $T$.

\begin{lemma} For any anisotropic torus $T$ over ${\Bbb Q}$, there are
field extensions $K_1,\cdots,$ $K_d$ of ${\Bbb Q}$, such that if
$S$ denotes the product of the norm 1 tori associated to $K_i$, then
there is an  embedding of $T$ into $S$.

\end{lemma} 

\noindent{\bf Proof :} As is well-known, there is an equivalence of
categories between tori over ${\Bbb Q}$ and 
finitely generated free ${\Bbb Z}$-module with an action of the
Galois group Gal$(\bar{\Bbb Q}/{\Bbb Q})$
 of the algebraic closure $\bar{\Bbb Q}$ of ${\Bbb Q}$. The
equivalence is given by associating to any torus $T$, 
its character group $X^*(T)$.
Choose a ${\Bbb Z}$-basis, say $\{ e_1,\cdots, e_d\}$ 
of the character group of   $T$. Suppose that $H_i$ is the subgroup of 
Gal$(\bar{\Bbb Q}/{\Bbb Q})$ which stabilises the vector $e_i$. The mapping
$g \rightarrow g\cdot e_i$ from Gal$(\bar{\Bbb Q}/{\Bbb Q})$ to
$X^*(T)$ gives a mapping from ${\Bbb Z}[G/H_i]$ to $X^*(T)$.
Summing over $i$,
we get a surjective map from $\sum_i {\Bbb Z}[G/H_i]$ to $X^*(T)$. This
gives an embedding from $T$ to $\prod_i R_{K_i/{\Bbb Q}}{\Bbb G}_m$ 
where $K_i$ is the fixed field of $H_i$. 
Since $T$ is anisotropic, the image of $T$ lands inside
the product of norm 1 tori.

\section{ Schanuel's conjecture implies conjecture 2} 

In this section we prove that
our conjecture 2 about closures in Euclidean topology 
of finitely generated subgroups in general tori
is a consequence of Schanuel's conjecture.

We should however add that although our approach in this paper is via 
Schanuel's conjecture, 
there is a possibility that there may be a simpler
proof for conjecture 2, just using the more primitive methods of 
Geometry of Numbers,
as we are dealing with a rather specific  number theoretic context.

We begin with the statement of Schanuel's conjecture which is one of the most
outstanding problems in transcendental number theory. In the statement of 
this conjecture as well as 
everywhere else in the paper, one means by $\log A $, for a complex number $A$ 
to be {\it any} complex number $B$ such that $\exp (B) = A$.

\begin{conj}(Schanuel's Conjecture) If $\alpha_1,\cdots,\alpha_n$ are algebraic numbers such that
$\log \alpha_1,\cdots,\log \alpha_n$ are linearly independent over 
${\Bbb Q}$, then $\log \alpha_1,\cdots,\log \alpha_n$ are algebraically 
independent over ${\Bbb Q}$.
\end{conj}

We will need the following lemma about number fields.

\begin{lemma}
Let $L$ be a number field which is Galois over ${\Bbb Q}$. Enumerate the
elements of the Galois group $G$ of $L$ over ${\Bbb Q}$ as 
$\sigma_1=1,\sigma_2,\cdots,\sigma_d$. For an
element $z$ of $L^*$, denote the various (not necessarily distinct) 
Galois conjugates of
$z$ by $z_1 =z, z_2=\sigma_2(z) ,\cdots, z_d = \sigma_d(z)$.
 Let $x$ be an element $ L^*$.  
Then there is a nonzero integer $m$ and a unit $\epsilon$ in (the ring of 
integers of) $L$,
such that whenever
$x_1^{n_1}x_2^{n_2}\cdots x_d^{n_d} $ is a unit in (the ring of integers of) 
$L$ for a 
$d$-tuple  $(n_1,n_2,\cdots,n_d)$ inside ${\Bbb Z}^d$,
$(x^m\epsilon)_1^{n_1}(x^m\epsilon)_2^{n_2}\cdots (x^m\epsilon)_d^{n_d} =1.$
The integer $m$ can be taken to be the order of the class group of $L$ times
the degree of $L$ over ${\Bbb Q}$, hence can be chosen to be independent of
$x$.
\end{lemma}

\noindent{\bf Proof :} Write the (fractional) ideal generated by $x$ as a 
product of prime ideals:
$$(x) = \wp_1^{m_1}\cdots \wp_r^{m_r}.$$
(We will assume that if a certain prime $\wp_i$ occurs in the above decomposition, so does 
any Galois conjugate of it, with exponent perhaps 0.)
Let $k$ be an integer such that each of the ideals $\wp_i^k$ is a 
principal ideal generated by, say $\varpi_i$.    
If $H_i$ denotes the subgroup of the Galois group which fixes the
prime ideal $\wp_i$, then the elements of $H_i$ will take $\varpi_i$ 
into itself up to a unit: $h(\varpi_i) = h_i\cdot \varpi_i$. 
Clearly $h \rightarrow h_i$ is a 1-cocycle on $H_i$ with values in the
group of units of $L^*$. Since $H_i$ is a finite group, a finite power
of the cocycle becomes a coboundary, i.e., there is a positive integer
$r$, and a unit $\nu_i$  such that 
$$h(\varpi_i^r) = h^r_i\cdot \varpi_i^r =  h(\nu_i)\nu_i^{-1} \varpi_i^r.$$ 
It follows that $\nu_i^{-1}\varpi_i^r$ is invariant under $H_i$. Thus 
we can choose generators $\pi_i$ of the principal ideals $\wp_i^{rk}$ 
in such a way that if 
an element of the Galois group takes one such ideal into another such,
then the same holds for the generators (and not just only up to units). 
From the equality of ideals,
$(x^{rk}) = (\pi_1^{m_1}\pi_2^{m_2}\cdots\pi_r^{m_r})$, 
there is a unit $\epsilon$ with, 
$x^{rk} \epsilon = \pi_1^{m_1}\pi_2^{m_2}\cdots\pi_r^{m_r}$. 
Now observe that if a product of certain elements of $L^*$ 
with any two {\it distinct} elements coprime (such as $\pi_i$'s) is a unit, 
then the product is in fact 1 (and is in some sense the `empty product').
From this it follows that   
$(x^{kr}\epsilon)_1^{n_1}(x^{kr}\epsilon)_2^{n_2} \cdots 
(x^{kr}\epsilon)_d^{n_d}=1.$ Finally, the proof given here works with the
choice of $m=rk$ to be the product of the order of the class group of $L$ 
and the degree of $L$ over ${\Bbb Q}$.

\vspace{4mm}
\begin{cor}{\bf (of the proof):} With the notation as in the lemma, given
elements $\{x^{(1)},x^{(2)},\cdots,x^{(n)}\}$ in $L^*$, 
there is an  integer $m$ and units
$\epsilon_i$ in the ring of integers of  $L$ 
such that for $y^{(i)}= (x^{(i)})^{m}\epsilon_i$, the
subgroup of $L^*$ generated by $\sigma_j(y^{(i)})$ does not contain any unit
of the ring of integers of $L$ other than 1.

\end{cor}

\begin{thm} Schanuel's conjecture implies conjecture 2.
\end{thm}

\noindent {\bf Proof :} As already noted in remark 2, it suffices
to prove conjecture 2 for cyclic subgroups. Furthermore, because of lemma 
3, we can assume that the anisotropic torus $S$ is the product of
the norm 1 tori: $S = \prod_{i=1}^m (R_{K_i/{\Bbb Q}}{\Bbb G}_m)^1$. 
We will further assume (by enlarging the anisotropic 
torus which does not affect the conclusion regarding closures) that the fields
$K_i$ are Galois over ${\Bbb Q}$, and by taking the compositum,
we assume that the fields $K_i$ are the same, say $L$, a Galois extension
of degree $d+1$ over ${\Bbb Q}$ which we will assume to be totally 
real. The case when $L$ has complex places is very similar, although
notationally more cumbersome.
 
For an element $x^{(i)} \in L$, we let
$x^{(i)}_j$, $j \in \{1,2,\cdots,d+1\}$,
denote the various Galois conjugates of $x^{(i)}$.

Let 
$x=(x^{(1)},x^{(2)},\cdots,x^{(m)}) \in S = 
\prod_{i=1}^m (R_{L/{\Bbb Q}}{\Bbb G}_m)^1$.   
Replacing $x^{(i)}$ by 
$y^{(i)} =(x^{(i)})^{m}\epsilon_i$ as in corollary 1, we can assume
that $y = (y^{(1)},y^{(2)},\cdots,y^{(m)})$ is such that the group generated by
the various Galois conjugates $y^{(i)}_j$ intersects the units in the ring of 
integers of $L$ in identity alone.

We note that a general algebraic character of $L^*$ is given by 
$z \rightarrow \prod_j(z_j)^{n_j}.$ 
Denote by $A$ the group of characters $\chi$ of $S=
\prod_{i=1}^m (R_{L/{\Bbb Q}}{\Bbb G}_m)^1$ 
such that $\chi(y)=1$. 
Let the rank of the abelian group $A$ be $dm-k$.
Therefore the
subgroup of $L^*$ generated by $\{y^{(i)}_j\}$ is a free abelian
group of rank $k$.

We will use the homomorphism with finite kernel:
$$
\begin{array}{ccc}
\left[ L \otimes {\Bbb R} \right]^1 & \rightarrow & {\Bbb R}^d \\
(x_1,\cdots,x_{d+1}) & \rightarrow & (\log |x_1|,\cdots,\log|x_{d}|)
\end{array}
$$
under which (by Dirichlet unit theorem) 
the group of units of the ring of integers of 
$L$ of norm 1, ${\cal O}_L^{*1}$ 
goes to a lattice $\Lambda$ in ${\Bbb R}^d$ with ${\Bbb R}^d/\Lambda$ 
compact.

Taking the direct sum of this homomorphism $m$ number of times, 
we get a homomorphism from $S({\Bbb R})$ to
${\Bbb R}^{dm}/\Lambda^m$, whose kernel is an arithmetic group
in $S({\Bbb R})$, to be denoted by $S({\Bbb Z})$. 
We will denote the image of  
$y = (y^{(1)},y^{(2)},\cdots,y^{(m)})$ in ${\Bbb R}^{dm}$ as 
$(\log (|y^{(1)}_1|), \cdots,\log (|y^{(1)}_d|),\cdots, \log (|y^{(m)}_1|), 
\cdots,
\log (|y^{(m)}_d|)).$ By lemma $1(c)$, to prove this theorem it suffices to
prove that the rank of the matrix 

$$
A = \left(
\begin{array}{ccccccc}
  a_{11} & \cdots & a_{1d}& \cdots   &a_{m1} & \cdots & a_{md} 
\\
  \ell(1)_{11} & \cdots & \ell(1)_{1d}& \cdots   &\ell(1)_{m1} & \cdots & 
\ell(1)_{md} \\
\cdot  & \cdots & \cdot  & \cdots & \cdot & \cdots &\cdot \\
 \cdot &  \cdots & \cdot & \cdots & \cdot &\cdots & \cdot\\
\cdot & \cdots & \cdot  & \cdots  & \cdot & \cdots & \cdot \\
  \ell(k-1)_{11} & \cdots & \ell(k-1)_{1d}& \cdots   &\ell(k-1)_{m1} & \cdots 
& \ell(k-1)_{md} 
\end{array}
\right)
 $$
is $k$ where  $ (a_{11}, \cdots , a_{1d}, \cdots   , a_{m1} , \cdots ,
 a_{md}) =
(r\log(|y^{(1)}_1|)+\ell_{11} , \cdots ,  r\log(|y^{(1)}_d|)+\ell_{1d} , 
\cdots   , r\log (|y^{(m)}_1|) +\ell_{m1} , \cdots ,
r\log (|y^{(m)}_d|))+\ell_{md})$,
$r$ is a nonzero integer, $(\ell_{11},\cdots, \ell_{1,d},
\cdots, \ell_{m1}, \cdots, \ell_{md})$ is a vector in $\Lambda^m$,
and where the 2nd to $k$th rows of this matrix represent 
$(k-1)$ ${\Bbb Z}$-linearly independent vectors in ${\Lambda^m}$. 

Since the rank of the matrix 

$$
B = \left(
\begin{array}{ccccccc}

  \ell(1)_{11} & \cdots & \ell(1)_{1d}& \cdots   &\ell(1)_{m1} & \cdots& 
\ell(1)_{md} \\
\cdot  & \cdots & \cdot  & \cdots & \cdot & & \cdots \\
 \cdot &  \cdots & \cdot & \cdots & \cdot & & \cdots \\
\cdot & \cdots & \cdot  & \cdots  & \cdot & & \cdots\\
  \ell(k-1)_{11} & \cdots & \ell(k-1)_{1d}& \cdots   &\ell(k-1)_{m1} 
& \cdots & \ell(k-1)_{md} 
\end{array}
\right)
 $$
is $(k-1)$, there is a $(k-1) \times (k-1)$ submatrix with 
nonzero determinant. After re-indexing the co-ordinates in ${\Bbb R}^{dm}$,
we assume that the first 
$(k-1) \times (k-1)$ submatrix of $B$ has rank $(k-1)$, i.e., has
nonzero determinant.

Since the rank of the group generated by $\{y^{(i)}_j\}$ is $k$, there is
at least one index, say $y^{(i_{0})}_{j_0}$, 
such that no power of it belongs to
the group generated by the $y's$ corresponding to the first $(k-1)$ entries
in the first row of $A$ (after re-indexing introduced above). 
Denote these $y's$ 
as $y_1,y_2,\cdots, y_{k-1}$, and the corresponding $\ell$'s as 
$\ell_1,\ell_2,\cdots, \ell_{k-1}$. Also, denote $y^{(i_{0})}_{j_0}$ as $y_k$. 

Let $C$ be the $k \times k$ submatrix of $A$ comprising of the
first $(k-1)$ columns of $A$, and the $k$-th column corresponding to
$y^{(i_{0})}_{j_0}$. We want to prove that $\det(C) \not = 0$.

Clearly $\det(C) = \sum_{i=1}^{k} [r\log(|y_i|) +\ell_i]\det (L_i)$ 
where $L_i$ is a $(k-1) \times (k-1)$ matrix consisting 
of log of units in $L$, and $\ell_i$ are also log of units in $L$.

It follows from Schanuel's conjecture and our hypothesis that
no (nonzero) power of $y_k$  belongs to the group generated by
$y_i, i=1,\cdots, k-1$, 
 that $\log(|y_k|)$
is algebraically independent over the subfield of ${\Bbb C}$ 
generated by log of algebraic units and the $\log (|y_i|), i=1,\cdots, k-1$. 

By our assumption,  the first 
$(k-1) \times (k-1)$ submatrix of $B$ has
nonzero determinant which is $\det(L_k)$, hence
$\det(C) = \sum [r\log(|y_i|)+\ell_i]\det(L_i)$
is non-zero (by algebraic independence of the $k$th term from the rest).

\section{Counter-example to a  more general question}

It is natural to ask  if  an analogue of conjecture 2 can be made   
more generally. 
The general question is about the algebraicity of the 
connected component of identity of the closure in Euclidean topology of a 
finitely generated subgroup of ${\Bbb C}^{*n}$ with algebraic 
co-ordinates where $ {\Bbb C}^{*n}$ is considered as the $2n$-dimensional
torus defined over ${\Bbb R}$ as the Weil restriction of scalars $R_{{\Bbb C}/
{\Bbb R}}{\Bbb G}_m^n$.
For example, can one drop the condition on the torus
$S$ in conjecture 2 being anisotropic over ${\Bbb Q}$, 
and instead of taking $S({\Bbb Z})$
which is a cocompact discrete subgroup of $S({\Bbb R})$ if $S$ is 
anisotropic, take any 
cocompact discrete subgroup $\Gamma$ of $S({\Bbb R})$ contained in 
$S({\Bbb Q})$?   A simple counter-example shows that this is not possible, 
shattering any hope for a simple answer to the general question above.

To construct the counter-example, take $S= {\Bbb G}^3_m$, the 3 dimensional 
split torus over  ${\Bbb Q}$. The principle behind the counter-example is
the well-known observation that although the determinant of a 
 skew-symmetric $n \times n$ matrix 
consisting of  logarithm of algebraic numbers is 0 if $n$ is odd,
the rows and columns could be linearly independent over ${\Bbb Q}$, 
such as for the matrix:
$$
 \left(
\begin{array}{ccc}
0 & \log 2 & \log 3 \\
-\log 2 & 0 & \log 5 \\
-\log 3 & -\log 5 & 0 
\end{array}
\right).
 $$
Since the determinant of the following matrix is nonzero,
$$
 \left(
\begin{array}{ccc}
0 & \log 2 & \log 3 \\
-\log 2 & 0 & \log 5 \\
\log 7 & 0 & \log 2 
\end{array}
\right),
 $$
it follows that the subgroup of ${\Bbb R}^{*3}$ generated by the 
elements, $x=(1,2,3), y=(1/2,1,5), z =(7, 1, 2)$, is a discrete
cocompact subgroup $\Gamma$ of ${\Bbb R}^{*3}$. However, 
the closure of the image of the cyclic group generated by
the image of $w= (3,5,1)$ in $\Gamma \backslash {\Bbb R}^{*3}$ is a 
2 dimensional topological torus (this follows from lemma 2(b)) 
which does not arise from an  algebraic subtorus of ${\Bbb R}^{*3}$ 
as is easy to see.

\vspace{4mm}

\noindent{\bf Remark: } 
The connected component of identity of the closure of a finitely generated 
subgroup of algebraic numbers of ${\Bbb C}^*$ is $\{1, {\Bbb R}^+, 
{\Bbb S}^1, {\Bbb C}^* \}$. We refer to theorem 1.10, page 56 of 
[W3] for a proof of
this assuming Schanuel's conjecture. (The subtlety lies in proving that 
algebraic points cannot be dense on a closed connected subgroup of
${\Bbb C}^*$ besides the above ones.)

\section{Non-abelian analogue} It seems very natural to extend the scope of 
the conjecture 2  by replacing the anisotropic 
torus $T$ by a general  algebraic group $G$ over 
${\Bbb Q}$. 

We recall that by a theorem due to  Borel and Harish-Chandra, 
for a reductive algebraic group $G$ over ${\Bbb Q}$ which is anisotropic over 
${\Bbb Q}$,
$G({\Bbb Z})\backslash G({\Bbb R})$ 
is compact. We would like to suggest the 
analogue of conjecture 2 to  assert that the
closure of the image 
in $G({\Bbb Z})\backslash G({\Bbb R})$ of a finitely generated subgroup 
$F$ of $G(\bar{{\Bbb Q}}_{\Bbb R})$ (where $\bar{{\Bbb Q}}_{\Bbb R}$ 
is the subfield of algebraic numbers in ${\Bbb R}$)
is of the form $\Gamma_H \backslash H({\Bbb R})$ for an algebraic subgroup $H$ of
$G$ defined over ${\Bbb Q}$ with
$\Gamma_H = H({\Bbb R}) \cap G({\Bbb Z})$. (Note that if the image of a 
subgroup
$H({\Bbb R})$ in  $G({\Bbb Z})\backslash G({\Bbb R})$ is closed, hence
compact, then $\Gamma_H = H({\Bbb R}) \cap G({\Bbb Z})$ is a 
co-compact lattice in $H({\Bbb R})$, and hence if $H$ is algebraic, it is 
defined over ${\Bbb Q}$
by the Borel density theorem.) Notice that we have not proved 
even for a torus (even after assuming Schanuel's conjecture), 
a theorem in this generality as we have always restricted ourselves to
finitely generated subgroups $F$ of the torus which are contained in the
group of ${\Bbb Q}$-rational points. This seems to have been necessary for
the proof of conjecture 2 given here.

We remark that our suggested analogue contains a consequence of
M. Ratner's theorem (the proof of the so-called Raghunathan conjecture)
as observed by Dani and Raghunathan, cf. Cor. 4.9 in [V] in a very special
case. It states that if a semi-simple group G over ${\Bbb R}$ (with 
$G({\Bbb R})$ non-compact) has two
distinct ${\Bbb Q}$ structures, with corresponding lattices
${\Gamma_1}$ and ${\Gamma_2}$, then $\Gamma_1 \cdot \Gamma_2$ is 
dense in $G({\Bbb R})$ (in the Euclidean topology).

\vspace{8mm}

\noindent{\bf Acknowledgement: }This note was conceived several years
back. It is  being published now in the hope of stimulating further research.
The author thanks B. Mazur, Gopal Prasad, T.N.Venkataramana  and  
M. Waldschmidt  for encouragement and helpful remarks.

}

\noindent{Harish-Chandra  Research Institute, 
Chhatnag Road, Jhusi, Allahabad 211019, 
INDIA}

\noindent{Email: dprasad@mri.ernet.in}

\end{document}